\documentclass[10pt]{article}
\ifx\pdfoutput\undefined
\usepackage[dvips]{graphicx}
\else
\usepackage[pdftex]{graphicx}
\fi
\catcode`\@=11
\newtheorem{theorem}{Theorem}

\newtheorem{proposition}{Proposition}

\newtheorem{corollary}{Corollary}
\def\demo{\noindent{\bf Proof .-}}
\def\section{\@startsection {section}{1}{\z@}{-3.5ex plus -1ex
minus-.2ex}{2.3ex plus .2ex}{\normalsize\bf}}

\begin{document}
\ifx\pdfoutput\undefined
\else
\DeclareGraphicsExtensions{.pdf,.gif,.jpg} % the formats we have images in
\fi
\begin{center}
{\Large\bf \textsc{The Stanley-Reisner ideals of polygons as set-theoretic complete intersections}}
\end{center}
\vskip.5truecm
\begin{center}
{Margherita Barile\footnote{Partially supported by the Italian Ministry of Education, University and Research.}\\ Dipartimento di Matematica,\\ Universit\`{a} di Bari ``Aldo Moro", Via E. Orabona 4,\\70125 Bari, Italy\footnote{e-mail: barile@dm.uniba.it}}
\end{center}
\begin{center}
{Naoki Terai\\ Department of Mathematics, Faculty of Culture and Education, Saga University,\\ Saga 840-8502, Japan}\footnote {e-mail: terai@cc.saga-u.ac.jp}
\end{center}
\vskip1truecm
\noindent
{\bf Abstract} We show that the Stanley-Reisner ideal of the one-dimensional simplicial complex whose diagram is an $n$-gon
 is always a set-theoretic complete intersection in any positive characteristic.  

\section*{Introduction and Preliminaries}
The {\it arithmetical rank} (ara) of an ideal $I$ in a commutative Noetherian ring $R$ is the minimal number $s$ of elements $a_1,\dots, a_s$ of $R$ such that $\sqrt I=\sqrt{(a_1,\dots, a_s)}$; one can express this equality by saying that $a_1,\dots, a_s$  {\it generate} $I$ {\it up to radical}. In general height\,$I\leq\,$ara\,$I$; if equality holds, $I$ is called a {\it set-theoretic complete intersection}. The ideal $I$ is called a {\it complete intersection} in the special case where $I$ is generated by height\,$I$ elements.
Let $X$ be a non-empty finite set of indeterminates over an algebraically closed field $K$. A {\it simplicial complex} on $X$ is a set $\Delta$ of subsets of $X$ such that for all $x\in X$, $\{x\}\in\Delta$ and whenever $F\in\Delta$ and $G\subset F$, then $G\in\Delta$. The elements of $\Delta$ are called  {\it faces}, whereas $X$ is called the {\it vertex set} of $\Delta$, and the elements of $X$ are called the {\it vertices} of $\Delta$. If $\Delta$ consists of all subsets of its vertex set, then it is called a {\it simplex}. The simplicial complex $\Delta$ can be associated with an ideal $I_{\Delta}$ of the polynomial ring $R=K[X]$, which is generated by all monomials whose support is not a face of $\Delta$; $I_{\Delta}$ is called the {\it  Stanley-Reisner ideal} of $\Delta$ (over $K$). Its minimal monomial generators are the products of the elements of the minimal non-faces of $\Delta$, and these are squarefree monomials. In fact, this construction provides a one-to-one correspondence between the simplicial complexes on $X$ and the squarefree monomial ideals of $K[X]$ that do not contain  elements of degree one. The quotient ring $K[\Delta]=K[X]/I_{\Delta}$ is called the {\it  Stanley-Reisner ring} of $\Delta$ (over $K$). Its Krull dimension is equal to $\max\{\vert F\vert\mid F\in\Delta\}$, and this number, lowered by one, is defined as the {\it dimension} of $\Delta$. \newline
If $\Delta$ is 1-dimensional, we can associate with $\Delta$ a graph $G(\Delta)$ on the same vertex set: its edges are the faces of $\Delta$ having exactly two elements.  
In this paper, we consider, for all integers $n\geq 3$, the simplicial complex $\Delta_n$ for which $G(\Delta_n)$  is  an $n$-gon. In this case $|X|=n$, say $X=\{x_1,\dots, x_n\}$, and 
$$I_n=I_{\Delta_n}=(x_ix_j\mid i\mbox{ is not adjacent to $j$ in the cycle }1,2,\dots, n).$$
\noindent
 Moreover, height\,$I_n=n-2$. We have that $I_3=(x_1x_2x_3)$ and $I_4=(x_1x_3, x_2x_4)$ are complete intersections.  In \cite{B1}, Example 1 it is shown that $I_5$ can be generated up to radical by 3 elements,  and in \cite{B0}, Example 2, that  $I_6$ can be generated up to radical by 4 elements. This shows that $I_n$ is a set-theoretic complete intersection for $n\in\{3,4,5,6\}$. In this paper we show that this property is true for all $n\geq3$ if the characterstic of $K$ is positive. We present a recursive procedure, which, starting from $n-3$ polynomials defining $I_{n-1}$ up to radical, constructs $n-2$ polynomials defining $I_n$ up to radical. To this end we develop a linear-algebraic technique, similar to the one in \cite{B0} and \cite{BT}, but this time we use a resultant instead of a determinant.
\par\smallskip\noindent
\section{The recursive construction}
Let char\,$K=p$. Fix an integer $N\geq 6$. In this section we show by induction on $n$, for $5\leq n\leq N$, that there exist $n-2$ polynomials $f_1, \dots, f_{n-2}$ defining $I_n$ up to radical and fulfilling certain conditions on their forms, which we will state below as (I), (II), (III). As a consequence, we will deduce that $I_n$ is a set-theoretic complete intersection for every $n\geq 5$, since we can take $N$ arbitrarily large. \newline
We first introduce some numerical invariants. For each $n=4,\dots, N$, fix a positive integer $r^{(n)}$ and, for each $n=5,\dots N$, fix an integer sequence 
$$p^{r^{(n)}}=\gamma_1^{(n)},\dots,\gamma_{n-3}^{(n)},$$
\noindent
in which $\gamma_1^{(n)}$ is the minimum, in such a way that, for all $n=5,\dots, N$, 
\begin{equation}\label{r} r^{(n-1)}>2r^{(n)}+N,\end{equation}
\noindent
and
\begin{equation}\label{gamma} \sum_{i=1}^{n-3}\gamma_i^{(n)}=p^{s^{(n)}},\end{equation}
\noindent
where
\begin{equation}\label{s} s^{(n)}=r^{(n)}+N.\end{equation}
Then set $\alpha^{(5)}=p^{r^{(5)}}$, $\delta^{(5)}=p^{r^{(4)}}$ and  $\beta^{(5)}= 1$, and pick  integers $\lambda^{(5)}>\delta^{(5)}$ and $\epsilon^{(5)}>\alpha^{(5)}\delta^{(5)}$. For all $n=6,\dots, N$, set
\begin{equation}\label{delta} \delta^{(n)}=p^{r^{(n-1)}},\end{equation}
and recursively define
\begin{equation}\label{alpha} \alpha^{(n)}=\alpha^{(n-1)}p^{s^{(n-1)}},\end{equation}
\begin{equation}\label{beta} \beta^{(n)}=\beta^{(n-1)}\lambda^{(n-1)},\end{equation}
where, for all $n=6,\dots, N$,
\begin{equation}\label{lambda} \lambda^{(n)}=\alpha^{(n-1)}.\end{equation}
Finally, for all $n=6,\dots, N$, set
\begin{equation}\label{epsilon} \epsilon^{(n)}=\alpha^{(n-1)}(\delta^{(n-1)}-p^{s^{(n-1)}}).\end{equation}
\noindent
Note that, for all $n=5,\dots, N$, in view of  (\ref{r}) and (\ref{delta}), we have
\begin{equation}\label{deltadelta} \delta^{(n-1)}=p^{r^{(n-2)}}>p^{2r^{(n-1)}+N}>p^{r^{(n-1)}}=\delta^{(n)}.\end{equation}
\noindent
Moreover, for all $n=6,\dots, N$,
\begin{equation}\label{lambdadelta}
\lambda^{(n)}=\alpha^{(n-1)}>p^{s^{(n-2)}}>p^{r^{(n-2)}}>p^{r^{(n-1)}}=\delta^{(n)},\end{equation}
\noindent
where we have used (\ref{lambda}), (\ref{alpha}), (\ref{s}), (\ref{r}) and (\ref{delta}).
\noindent
From (\ref{epsilon}), (\ref{delta}),  (\ref{alpha}), (\ref{s}) and (\ref{r})   we deduce that, for all $n=5, \dots, N$
\begin{eqnarray}\label{epsilonpositive} 
\epsilon^{(n)}&=&\alpha^{(n-1)}(p^{r^{(n-2)}} -p^{s^{(n-1)}})\nonumber\\
&=&\alpha^{(n-1)}p^{s^{(n-1)}}(p^{r^{(n-2)}-s^{(n-1)}} -1)\nonumber\\
&=&\alpha^{(n)}(p^{r^{(n-2)}-r^{(n-1)}-N}-1)\nonumber\\
&>&\alpha^{(n)}(p^{2r^{(n-1)}+N-r^{(n-1)}-N}-1)\nonumber\\
&=&\alpha^{(n)}(p^{\delta^{(n)}}-1)>\alpha^{(n)}\delta^{(n)}.
\end{eqnarray}
\noindent
In the sequel, for the sake of simplicity, we will throughout omit the superscript $(n-1)$, so that, e.g., $\alpha$ will stand for $\alpha^{(n-1)}$. 
Fix an integer $n$ such that $6\leq n\leq N$. We suppose that $I_{n-1}$ can be generated, up to radical, by polynomials $g_1, \dots, g_{n-3}\in K[x_1,\dots,x_{n-1}]$ such that 
\begin{eqnarray}
g_1&=&a_{1,2}x_2^{\lambda\beta}+\sum_{j=3}^{n-3}a_{1,j}x_j^{\lambda},\label{g1}\\
g_2&=&a_{2,2}x_2^{\beta}+\sum_{j=3}^{n-3}a_{2,j}x_j\label{g2},
\end{eqnarray}
\noindent
 and, for all $i=3,\dots, n-4$, 
\begin{equation}\label{gi}g_i=\sum_{j=3}^{n-3}a_{i,j}x_j,\end{equation}
\noindent
and, finally
\begin{equation}\label{gn3}g_{n-3}=\sum_{j=3}^{n-3}a_{n-3,j}x_j^{\alpha}+x_1x_{n-2}^{\alpha},\end{equation}
 where, for all indices $i$ and $j$,  $a_{i,j}\in K[x_1,\dots,x_{n-1}]$. Furthermore we assume that   
\begin{list}{}{}
\item{(I)} for all $i=1,\dots, n-4$,
$$a_{i,i+1}=x_{n-1}^{\gamma_i}+b_i,$$
\noindent
where  
$$b_i\in (x_3^{\delta},\dots, x_{n-2}^{\delta});$$
\noindent
\item{(II)} for all indices $i=2,\dots, n-4$ and $j=3,\dots, n-3$ such that  $j\neq i+1$,
 we have
$$a_{i,j}\in(x_3^{\delta},\dots, x_{n-2}^{\delta}),$$
moreover
\begin{equation}\label{d} a_{2,2}\in(x_3^{\delta},\dots, x_{n-2}^{\delta}),\quad\mbox{say }a_{2,2}=\sum_{j=3}^{n-2}d_jx_j^{\delta}.\end{equation}
\item{(III)} for all indices $j=3,\dots, n-3$,
$$a_{n-3,j}\in I_{n-1}\cap(x_3^{\epsilon},\dots,  x_{n-3}^{\epsilon}).
$$
\end{list}
\smallskip
\noindent
We will construct $n-2$ polynomials $f_1,\dots, f_{n-2}\in K[x_1,\dots, x_n]$ that have the same form as $g_1,\dots, g_{n-3}$ and generate $I_n$ up to radical.\newline 
First we rewrite $g_1$ and $g_2$. Equality (\ref{g1new}) follows from  (\ref{g1}), whereas (\ref{g2new}) is a consequence of (\ref{g2}) and (\ref{d}).
\begin{eqnarray}
g_1&=&a_{1,2}x_2^{\lambda\beta}+\sum_{j=3}^{n-3}a'_{1,j}x_j,\label{g1new}\\
g_2&=&\sum_{j=3}^{n-2}a'_{2,j}x_j,\label{g2new}
\end{eqnarray}
\noindent
where, for all indices $j=3,\dots, n-3$,
\begin{equation}\label{a1prime} a'_{1,j}=a_{1,j}x_j^{\lambda-1}\in (x_3^{\delta}, \dots, x_{n-3}^{\delta}),\end{equation}
\noindent
(which is true because from the definition of $\lambda^{(5)}$ and from (\ref{lambdadelta}) we know that $\lambda >\delta$), and, for $j=4,\dots, n-2$,
\begin{equation}\label{a2prime}a'_{2,j}=a_{2,j}+d_jx_2^{\beta}x_j^{\delta}\in (x_3^{\delta}, \dots, x_{n-2}^{\delta}),\end{equation}
\noindent
which is true by virtue of (II). Here, according to (\ref{g2}), we have set $a_{2,n-2}=0$. Moreover, by (I), 
\begin{equation}\label{a23prime} a'_{2,3}=a_{2,3}+d_3x_2^{\beta}x_3^{\delta}=x_{n-1}^{\gamma_2}+b'_2,\end{equation}
\noindent
where we have set
$$b'_2=b_2+d_3x_2^{\beta}x_3^{\delta}\in (x_3^{\delta}, \dots, x_{n-2}^{\delta}).$$
\noindent 
For all $i=1,\dots, n-3$, let $\tilde g_i\in K[x_1,\dots, x_{n-1}, y_2,\dots, y_{n-2}]$ be defined as follows. 
Let 
\begin{eqnarray}
\tilde g_1&=&a_{1,2}y_2+\sum_{j=3}^{n-3}a'_{1,j}y_j,\\\label{tildeg1}
\tilde g_2&=&\sum_{j=3}^{n-2}a'_{2,j}y_j.\label{tildeg2}
\end{eqnarray}
Then, for $i=3,\dots, n-4$, set
\begin{equation}\label{tildegi}\tilde g_i=\sum_{j=3}^{n-3}a_{i,j}y_j,\end{equation}
\noindent
and set
\begin{equation}\label{tildegn3}\tilde g_{n-3}=\sum_{j=3}^{n-3}a_{n-3,j}y_j^{\alpha}+x_1y_{n-2}^{\alpha}.\end{equation}
\noindent
In the sequel the apex will replace the superscript $(n)$, so that, e.g., $\alpha'$ will stand for $\alpha^{(n)}$.  
We now define $f_1,\dots, f_{n-3}\in K[x_1,\dots, x_n]$ by setting
\begin{equation}\label{f1}f_1=g_{n-3}+x_n^{\gamma'_1}x_2^{\alpha\lambda\beta},\end{equation}
\noindent
and, for all $i=2,\dots, n-3$,
\begin{equation}\label{fi}f_i=g_{i-1}+x_n^{\gamma'_i}x_{i+1}.\end{equation}
 Moreover, we set
$$\tilde f_1= \tilde g_{n-3}+x_n^{\gamma'_1}y_2^{\alpha},$$
\noindent
and, for all $i=2,\dots, n-3$,
$$\tilde f_i= \tilde g_{i-1}+x_n^{\gamma'_i}y_{i+1}.$$
\noindent
Then, in view of (22), (\ref{tildeg2}) and (\ref{tildegi}), we have that
$\tilde f_2,\dots, \tilde f_{n-3}$ are linear polynomials in $y_2,\dots, y_{n-2}$ and, by (\ref{tildegn3}),
\begin{equation}\label{tildef1}
\tilde f_1=x_n^{\gamma'_1}y_2^{\alpha}+\sum_{j=3}^{n-3}a_{n-3,j}y_j^{\alpha}+x_1y_{n-2}^{\alpha}\end{equation}
\noindent
is a homogeneous polynomial in $y_2,\dots, y_{n-2}$. By (\ref{fi}) we also have, in view of (22),
$$\tilde f_2=\sum_{j=2}^{n-2}c_{2,j}y_j,$$
\noindent
where
$$c_{2,j}=\left\{
\begin{array}{ll} 
a_{1,2}&\mbox{if }j=2,\\
a'_{1,3}+x_n^{\gamma'_2}&\mbox{if }j= 3,\\
a'_{1,j}&\mbox{if }4\leq j\leq n-3,\\
0&\mbox{if }j= n-2,
\end{array}
\right.
$$
\noindent
and, in view of (\ref{tildeg2}),
$$\tilde f_3=\sum_{j=3}^{n-2}c_{3,j}y_j,$$
\noindent
where
$$c_{3,j}=\left\{
\begin{array}{ll} a'_{2,j}&\mbox{if }j\neq 4\\
a'_{2,4}+x_n^{\gamma'_3}&\mbox{if }j= 4.
\end{array}
\right.
$$
\noindent
Finally, for $i=4,\dots, n-3$ we have
$$\tilde f_i=\sum_{j=2}^{n-2}c_{i,j}y_j,$$
\noindent
where
$$c_{i,j}=\left\{
\begin{array}{ll} a_{i-1,j}&\mbox{if }j\neq i+1, n-2,\\
a_{i-1,i+1}+x_n^{\gamma'_i}&\mbox{if }j= i+1,\\
0&\mbox{if }j= n-2.
\end{array}
\right.
$$
\noindent
Here we have set $a_{n-4,n-2}=0$. 
Let $S$ denote the resultant of $\tilde f_1,\dots,\tilde f_{n-3}$ as polynomials in the indeterminates $y_1,\dots, y_{n-2}$. Then, by \cite{J}, Proposition 5.4.4, 
\begin{equation}\label{R}S=\tilde f_1(\Delta_1,\dots, \Delta_{n-3}),\end{equation}
\noindent
where $\Delta_1,\dots, \Delta_{n-3}$ are such that $\det B=\sum_{i=1}^{n-3}\Delta_iT_i$; here $T_1,\dots, T_{n-3}$ are indeterminates over $K[x_1,\dots, x_n]$ and 
\begin{eqnarray*}&&B=\left(\begin{array}{cccc}
c_{2,2}&\cdots& c_{2,n-3}&c_{2,n-2}\\
\vdots&&\vdots&\vdots\\
c_{n-3,2}&\cdots& c_{n-3,n-3}&c_{n-3,n-2}\\
T_1&\cdots&T_{n-4}&T_{n-3}
\end{array}
\right)=\\\vphantom{X}\\
&&\!\!\!\!\!\!\!\!\!\!\!
\!\!\!\!\!\!\!\!\!\!\!
\!\!\!\!\!\!\!\!\!\!\!
\!\!\!\!\!\!\!\!\!\!\!
\!\!\!\!\!\!\!\!\!\!\!\left(\begin{array}{ccccccc}
x_{n-1}^{\gamma_1}+b_1& a'_{1,3}+x_n^{\gamma'_2}& a'_{1,4}&\cdots&a'_{1,n-4}&a'_{1,n-3}&0\\\ \\
0& x_{n-1}^{\gamma_2}+b'_2& a'_{2,4}+x_n^{\gamma'_3}&\cdots&a'_{2,n-4}&a'_{2,n-3}&a'_{2,n-2}\\\ \\
\vdots&\vdots&\ddots&\ddots&&&\vdots\\
0& a_{n-6,3}&\cdots&x_{n-1}^{\gamma_{n-6}}+b_{n-6} &a_{n-6,n-4}+x_n^{\gamma'_{n-5}}&a_{n-6,n-3}&0\\\ \\
0& a_{n-5,3}&a_{n-5,4} &\cdots&x_{n-1}^{\gamma_{n-5}}+b_{n-5}&a_{n-5,n-3}+x_n^{\gamma'_{n-4}}&0\\\ \\
0& a_{n-4,3}& a_{n-4,4}&\cdots&a_{n-4,n-4}&x_{n-1}^{\gamma_{n-4}}+b_{n-4}&x_n^{\gamma'_{n-3}}\\\ \\
T_1&T_2&T_3&\cdots&\cdots&T_{n-4}&T_{n-3}
\end{array}
\right).\\
\end{eqnarray*}
\noindent
In the second equality we have applied condition (I) and definition (\ref{a23prime}). 
Note that, by (I) and (\ref{gamma}),
\begin{eqnarray*}
\Delta_1&=&\bar\Delta_1+x_n^{p^{s'}-\gamma'_1},\\
\Delta_{n-3}&=&\bar\Delta_{n-3}+x_{n-1}^{p^s},
\end{eqnarray*}
\noindent
where $\bar\Delta_1,\bar\Delta_{n-3}\in (x_3^{\delta},\dots,  x_{n-2}^{\delta})$: this follows from conditions (I) and (II), (\ref{a1prime}), (\ref{a2prime}), since $\bar\Delta_1$ and $\bar\Delta_{n-3}$ belong to the ideal generated by the elements $b_1,b_2', b_3,\dots, b_{n-4}$, the elements  $a'_{1,j}$ such that $3\leq j\leq n-3$, the elements $a'_{2,j}$ such that $4\leq j\leq n-2$, and the elements $a_{i,j}$ such that $3\leq i\leq n-4$, $3\leq j\leq n-3$, and $j\neq i+1$. Therefore 
\begin{equation}\label{x1Delta}
x_n\bar\Delta_1^{\alpha}\in I_n\cap(x_3^{\alpha\delta},\dots,  x_{n-2}^{\alpha\delta})\end{equation}
\noindent
 and 
\begin{equation}\label{x3Delta}
x_1\bar\Delta_{n-3}^{\alpha}\in  I_n\cap(x_3^{\alpha\delta},\dots,  x_{n-2}^{\alpha\delta}).
\end{equation}
\noindent
 Now, in view of  (\ref{tildef1}) and (\ref{R}) 
\begin{eqnarray*}
S&=&x_n^{\gamma'_1}\Delta_1^{\alpha}+\sum_{j=3}^{n-3}a_{n-3,j}\Delta_{j-1}^{\alpha}+
x_1\Delta_{n-3}^{\alpha}\\&=&
F+x_n^{\gamma'_1(1-\alpha)+\alpha p^{s'}}+x_1x_{n-1}^{\alpha p^s},
\end{eqnarray*}
where 
$$F=x_n^{\gamma'_1}\bar\Delta_1^{\alpha}+\sum_{j=3}^{n-3}a_{n-3,j}\Delta_{j-1}^{\alpha}+
x_1\bar\Delta_{n-3}^{\alpha}.$$
Now, in view of (III) and (\ref{epsilonpositive}),
$$\sum_{j=3}^{n-3}a_{n-3,j}\Delta_{j-1}^{\alpha}\in I_{n-1}\cap(x_3^{\epsilon},\dots,  x_{n-3}^{\epsilon})\subset I_n\cap(x_3^{\alpha\delta},\dots,  x_{n-2}^{\alpha\delta}).$$
\noindent
From this relation, (\ref{x1Delta}) and (\ref{x3Delta}) we finally conclude that
\begin{equation}\label{4'}
F\in I_n\cap (x_3^{\alpha\delta},\dots,  x_{n-2}^{\alpha\delta}).
\end{equation}
\noindent
 Set
\begin{equation}\label{fn2}f_{n-2}=S-x_n^{\gamma'_1(1-\alpha)+\alpha p^{s'}}=F+x_1x_{n-1}^{\alpha p^s}.\end{equation}\noindent
\begin{proposition}\label{prop1} If the polynomials $g_1,\dots, g_{n-3}$ generate $I_{n-1}$ up to radical, then the polynomials $f_1,\dots, f_{n-2}$ generate $I_n$ up to radical. 
\end{proposition}
\demo We have to show that $I_n=\sqrt{(f_1,\dots, f_{n-2})}$. The inclusion $\supset$ is true because $f_1,\dots, f_{n-2}\in I_n$: the latter statement follows from the definitions of the polynomials $f_i$ and the fact that 
\begin{equation}\label{In}I_n=I_{n-1}K[x_1,\dots, x_n]+(x_1x_{n-1}, x_2x_n, \dots, x_{n-2}x_n).\end{equation}
\noindent
We prove the inclusion $\subset$. By Hilbert's Nullstellensatz, it is sufficient to show that whenever ${\bf x}=(x_1,\dots, x_n)\in K^n$ annihilates all $f_i$, then it annihilates all elements of $I_n$.  In the rest of the proof, for the sake of simplicity, in our notation we will identify each polynomial with its value at $\bf x$. Assume that $f_i=0$ for all $i=1,\dots, n-2$. We distinguish two cases. First suppose that $S=0$. Then, in view of (\ref{fn2})  from $f_{n-2}=0$ we deduce that $x_n=0$. But then, in view of (\ref{f1}) and (\ref{fi}), from $f_1=\cdots =f_{n-3}=0$ we derive that  $g_1=\cdots=g_{n-3}=0$. This implies that all polynomials in $I_{n-1}$ vanish at $\bf x$. Note that (\ref{4'}) and (\ref{In}) imply
\begin{equation}\label{5'}
F\in I_{n-1}K[x_1,\dots, x_n]+(x_n).
\end{equation}
\noindent
Hence $F=0$, so that, in view of (\ref{fn2}), from $f_{n-2}=0$ we conclude that $x_1x_{n-1}=0$. In view of (\ref{In}), this shows that all elements of $I_n$ vanish at $\bf x$. Now suppose that $S\ne0$. Then, by \cite{VDW}, p.~15, from 
$f_1=\cdots =f_{n-3}=0$ we conclude that $x_2=\cdots=x_{n-2}=0$. But $I_{n-1}\subset (x_2,\dots, x_{n-2})$, so that all elements of $I_{n-1}$ vanish at $\bf x$. Moreover, in view of (\ref{4'}), we have that $F=0$. Again from $f_{n-2}=0$ we deduce that $x_1x_{n-1}=0$. In view of (\ref{In}), we conclude that  all elements of $I_n$ vanish at $\bf x$. This completes the proof.
\par\smallskip\noindent
 Proposition \ref{prop1} provides the recursive step of the construction by which we intend to prove that $I_n$ is a set-theoretic complete intersection. In order to achieve this goal, we still have to 
\begin{list}{}{}
\item{(a)} provide three polynomials $g_1, g_2, g_3\in K[x_1, \dots, x_5]$ that generate $I_5$ up to radical and ensure that they fulfil conditions (I), (II) and (III), and
\item{(b)} show that the polynomials $f_1, \dots f_{n-2}$ fulfil conditions (I), (II) and (III), as well. 
\end{list}
\noindent
Task (a) is accomplished if we take $\gamma_1=\gamma_1^{(5)}$, $\gamma_2=\gamma_2^{(5)}$, $\delta=\delta^{(5)}$ and $\alpha=\alpha^{(5)}$, and we consider the polynomials 
\begin{eqnarray*}
g_1&=&x_5^{\gamma_1}x_2+x_1x_3^{\lambda},\\
g_2&=&x_4^{\lambda}x_2+x_5^{\gamma_2}x_3,\\
g_3&=&x_1x_4^{\alpha}.
\end{eqnarray*}
 These polynomials arise from a modification of those presented in \cite{B1}, Example 1. The proof can be easily performed by applying Hilbert's Nullstellensatz and considering the cases where $x_1=0$ and $x_4=0$ separately. Since $a_{1,2}=x_5^{\gamma_1}$, $a_{2,3}=x_5^{\gamma_2}$, condition (I) holds; furthermore, since $a_{2,2}=x_4^{\lambda}$, condition (II) holds; finally, since $a_{3,3}=0$, condition (III) holds, too.  We now handle (b). To this end, we first need to write $f_1,\dots, f_{n-2}$ in a suitable way.\newline
By (\ref{gn3}) and (\ref{f1}), we have
\begin{equation}\label{f1'}f_1=\tilde a_{i,2}x_2^{\alpha\lambda\beta}+\sum_{j=3}^{n-2}\tilde a_{1,j}x_j^{\alpha},\end{equation}
\noindent 
where 
\begin{equation}\label{star}\tilde a_{1,j}=\left\{
\begin{array}{ll} 
x_n^{\gamma'_1}&\mbox{if }j=2,\\
a_{n-3,j}&\mbox{if }3\leq j\leq n-3,\\
x_1&\mbox{if }j=n-2.
\end{array}
\right.
\end{equation}
\noindent
Moreover, by (\ref{g1}) and (\ref{fi}),
\begin{equation}\label{f2}f_2=\tilde a_{2,2}x_2^{\lambda\beta}+\sum_{j=3}^{n-2}\tilde a_{2,j}x_j,\end{equation}
\noindent 
where, by (I),
\begin{equation}\label{star2}\tilde a_{2,j}=\left\{
\begin{array}{ll} 
a_{1,2}=x_{n-1}^{\gamma_1}+b_1&\mbox{if }j=2,\\
x_n^{\gamma'_2}+a_{1,3}x_3^{\lambda-1}&\mbox{if }j=3,\\
a_{1,j}x_j^{\lambda-1}&\mbox{if }4\leq j\leq n-3,\\
0&\mbox{if }j=n-2,
\end{array}
\right.
\end{equation}
\noindent
and, by (\ref{fi}), (\ref{g2new}), (\ref{a2prime}) and (\ref{a23prime}), 
\begin{equation}\label{f3}f_3=\sum_{j=3}^{n-2}\tilde a_{3,j}x_j,\end{equation}
\noindent 
where
\begin{equation}\label{star3}\tilde a_{3,j}=\left\{
\begin{array}{ll} 
a_{2,j}+d_jx_2^{\beta}x_j^{\delta}&\mbox{if }j\neq 4, n-2\\
x_n^{\gamma'_3}+a_{2,4}+d_4x_2^{\beta}x_4^{\delta}&\mbox{if }j=4,\\
d_{n-2}x_2^{\beta}x_{n-2}^{\delta}&\mbox{if }j=n-2.
\end{array}
\right.
\end{equation}
\noindent
For $i=4,\dots, n-3$, we further have
\begin{equation}\label{fi'}
f_i=\sum_{j=3}^{n-2}\tilde a_{i,j}x_j,
\end{equation}
\noindent
where, by (\ref{gi}), and (\ref{fi}),
\begin{equation}\label{tilde}\tilde a_{i,j}=\left\{
\begin{array}{ll} a_{i-1,j}&\mbox{if }j\neq i+1,n-2,\nonumber\\
x_n^{\gamma'_i}+a_{i-1, i+1}&\mbox{if }j= i+1,\\
0&\mbox{if }j= n-2, i\ne n-3.
\end{array}
\right.
\end{equation}
\noindent 
Here, according to (\ref{gi}), we have set $a_{n-4,n-2}=0$, so that $\tilde a_{n-3,n-2}=x_n^{\gamma'_{n-3}}$. 
Finally, by  (\ref{4'}) and (\ref{fn2}) we have that 
\begin{equation}\label{fn2'}
f_{n-2}=\sum_{j=3}^{n-2}a_{n-2,j}x_j^{\alpha p^s}+x_1x_{n-1}^{\alpha p^s},
\end{equation}
\noindent
where, for all $j=3,\dots, n-2$,
\begin{equation}\label{doublestar}
a_{n-2,j}\in  I_n\cap (x_3^{\epsilon'},\dots,  x_{n-2}^{\epsilon'}),
\end{equation} 
and, according to (\ref{epsilon}), $\epsilon'=\alpha\delta-\alpha p^s$. Comparing (\ref{g1}) with (\ref{f1'}), (\ref{g2}) with (\ref{f2}), (\ref{gi}) with (\ref{f3}) and (\ref{fi'}), (\ref{gn3}) with (\ref{fn2'}), we see that $f_1,\dots, f_{n-2}$ are polynomials of the same form as $g_1,\dots, g_{n-3}$ with $n-2$ instead of $n-3$, $\lambda'=\alpha$ instead of $\lambda$ (see (\ref{lambda})), $\beta'=\lambda\beta$ instead of $\beta$ (see (\ref{beta})) and $\alpha'=\alpha p^s$  instead of $\alpha$  (see (\ref{alpha})). We show that conditions (I), (II) and (III) are fulfilled by $f_1,\dots, f_{n-2}$ with respect to these new data.   From (\ref{deltadelta}) we know that $\delta>\delta'$. Hence, in view of (II), for $i=4,\dots, n-4$, we have
$$a_{i-1,i+1}\in(x_3^{\delta}, \dots, x_{n-2}^{\delta})\subset (x_3^{\delta'}, \dots, x_{n-1}^{\delta'}).$$
Therefore, in view of (\ref{tilde}), the coefficients $\tilde a_{i,i+1}$, for $i=4,\dots, n-4$, fulfil condition (I) with $n$ instead of $n-1$, $\tilde b_i= a_{i-1,i+1}$ instead of $b_i$, $\gamma'_i$ instead of $\gamma_i$ and $\delta'$ instead of $\delta$. This is also true for $i=n-3$, with $\tilde b_{n-3}=0$, since $\tilde a_{n-3,n-2}=x_n^{\gamma'_{n-3}}$, and for $i=1$,  with $\tilde b_1=0$, since by (\ref{star}) $\tilde a_{1,2}=x_n^{\gamma'_1}$.   It is also true for $i=2$, with $\tilde b_2=a_{1,3}x_3^{\lambda-1}$, since, by virtue of (\ref{star2}), $\tilde a_{2,3}=x_n^{\gamma'_2}+a_{1,3}x_3^{\lambda-1}$ and, by (\ref{lambdadelta}) and (\ref{deltadelta}), $\lambda-1>\delta'$.  From (\ref{star3}) we have
$\tilde a_{3,4}=x_n^{\gamma'_3}+a_{2,4}+d_4x_2^{\beta}x_4^{\delta},$
where $a_{2,4}\in (x_3^{\delta}, \dots, x_{n-2}^{\delta})$: this follows from (II) if $n>6$, on the other hand, for $n=6$, we have that $a_{2,4}=a_{n-4,n-2}=0$. Since, by (\ref{deltadelta}), $\delta>\delta'$, it follows that condition (I) is also true for $i=3$ with $\tilde b_3=
 a_{2,4}+d_4x_2^{\beta}x_4^{\delta}$. 
 This establishes condition (I). 
 Next we show that condition (II) holds for the coefficients $\tilde a_{i,j}$ with $\delta'$ instead of $\delta$. 
From (\ref{tilde}) and (II) we see that 
\begin{equation}\label{square}\tilde a_{i,j}=a_{i-1, j}\in(x_3^{\delta},\dots, x_{n-2}^{\delta})\subset(x_3^{\delta'},\dots, x_{n-1}^{\delta'})
\end{equation}
 for all $i=4,\dots, n-3$ and $j=3,\dots, n-3$ such that $j\neq i$ and $j\ne i+1$. This also holds for $i=4,\dots, n-4$, $j=n-2$, since $\tilde a_{i,n-2}=0$. 
  From (II) and (\ref{star3}), we have
$$\tilde a_{3,j}=a_{2,j}+d_jx_2^{\beta}x_j^{\delta}\in (x_3^{\delta}, \dots, x_{n-2}^{\delta})\subset (x_3^{\delta'}, \dots, x_{n-1}^{\delta'})$$
\noindent
for $j=5,\dots, n-2$, where we have set $a_{2,n-2}=0$. Moreover, since $\lambda>\delta'$,
from (\ref{star2}) we derive that 
$$\tilde a_{2,j}=a_{1,j}x_j^{\lambda-1}\in (x_4^{\lambda-1}, \dots, x_{n-3}^{\lambda-1})\subset(x_3^{\delta'}, \dots, x_{n-1}^{\delta'})$$
\noindent
for all $j=4,\dots, n-3$, whereas $\tilde a_{2,n-2}=0.$  This establishes condition (II) for $i\neq j$. 
Now, from (I) and (\ref{tilde}) we deduce that, for all $i=4,\dots, n-3$, 
\begin{eqnarray*}
\tilde a_{i,i}&=&a_{i-1,i}=x_{n-1}^{\gamma_{i-1}}+b_{i-1}\\
&\in&(x_3^{\delta}, \dots, x_{n-2}^{\delta}, x_{n-1}^{\gamma_{i-1}})\subset(x_3^{\delta'}, \dots, x_{n-2}^{\delta'}, x_{n-1}^{\delta'}),
\end{eqnarray*}
\noindent
where we have used the fact that $\gamma_{i-1}\geq \gamma_1=\delta'$. 
Similarly, from (I) and (\ref{star2}) we deduce that
$$
\tilde a_{2,2}=x_{n-1}^{\gamma_1}+b_1
\in(x_3^{\delta'}, \dots, x_{n-2}^{\delta'}, x_{n-1}^{\delta'}),
$$
\noindent
and from (I) and (\ref{star3}), since $\gamma_2\geq \gamma_1$,
$$
\tilde a_{3,3}=a_{2,3}+d_3x_2^{\beta}x_3^{\delta}=x_{n-1}^{\gamma_2}+d_3x_2^{\beta}x_3^{\delta}+b_2
\in(x_3^{\delta'}, \dots, x_{n-2}^{\delta'}, x_{n-1}^{\delta'}).
$$
\noindent
This establishes condition (II).  Finally, note that (\ref{doublestar})  implies that the coefficients $a_{n-2,j}$ fulfil condition (III) with 
$\epsilon'$ instead of $\epsilon$.  We have just proven the following result.
\begin{theorem}\label{final} Suppose that char\,$K=p>0$. Then, for all $n\geq 5$, the ideal $I_n$ of $K[x_1,\dots, x_n]$ is a set-theoretic complete intersection.
\end{theorem}
\section{Some consequences}
The simplicial complexes $\Delta_n$ occur in some classification theorems, together with the simplicial complexes $\Lambda_n$ for which $G(\Lambda_n)$ is the straight path $x_1, \dots, x_n$. From \cite{S}, Theorem 5.1, we know that  the Stanley-Reisner ring $K[\Delta]$ of a 1-dimensional simplicial complex $\Delta$ is Gorenstein (i.e., $\Delta$ is a Gorenstein complex over $K$) if and only if $\Delta=\Delta_n$ for some $n\geq 3$, or $\Delta=\Lambda_2$, or $\Delta=\Lambda_3$. The first case is the one where the $a$-invariant of $K[\Delta]$ is zero. The Stanley-Reisner ideal $I_{\Lambda_2}$ is the zero ideal of $K[x_1,x_2]$, and $I_{\Lambda_3}= (x_1x_3)\subset K[x_1,x_2, x_3]$: both ideals are principal. As a consequence of Theorem \ref{final} we thus have:
\begin{corollary} Suppose that char\,$K>0$. If $\Delta$ is a 1-dimensional Gorenstein complex over $K$, then $I_{\Delta}$ is a set-theoretic complete intersection. 
\end{corollary}
More in general, $I_{\Lambda_n}$ is a set-theoretic complete intersection: this can be derived, e.g., from \cite{BT}, Corollary 2, by a trivial inductive argument.
In \cite{TY} Terai and Yoshida call $\Delta$ a {\it locally complete intersection} complex if, for the link of every vertex, the Stanley-Reisner ideal is  a complete intersection (recall that the link of a vertex $x$ is the subcomplex induced on the set of vertices distinct from $x$ that lie in the same face as $x$). In \cite{TY}, Theorem 1.15, they show that a non-empty simplicial complex $\Delta$ is a locally complete intersection if and only if $G(\Delta)$ is a finite disjoint union of ``polygons'', ``straight paths'', ``points'' and  complete intersection complexes of dimension at least 2, and that under this assumption, whenever dim\,$\Delta= 1$,  $\Delta$ is Cohen-Macaulay over $K$ if and only if $G(\Delta)$ is connected. It is well known that a complete intersection is Cohen-Macaulay. On the other hand, if it is disconnected,  it is well known that depth\,$ K[\Delta]=1$.
As a consequence of the Auslander-Buchsbaum formula (see, e.g., \cite{BH}, Theorem 1.3.3), we thus have projdim\,$K[\Delta]=n-1$,
where we recall that $n=\vert X\vert$. Moreover, projdim\,$K[\Delta] \leq\,$ ara\,$ I_{\Delta}$ by \cite{L}, and
ara\,$ I_{\Delta} \le n-1$ by \cite{EE}, Theorem 2.
In view of Theorem \ref{final} we thus obtain:
\begin{corollary} Suppose that char\,$K>0$. If $\Delta$ is a Cohen-Macaulay locally complete intersection complex over $K$, then $I_{\Delta}$ is a set-theoretic complete intersection. More in general, if $\Delta$ is a locally complete intersection complex over $K$, then 
$${\rm projdim}\, K[\Delta] =\, {\rm ara}\, I_{\Delta} = n-1.$$ 
\end{corollary}

\end{document}